\documentclass[11pt]{amsart}
\usepackage{amsfonts}
\usepackage{amsmath,amsthm,latexsym,amssymb,graphicx}
\usepackage{graphicx,color,enumerate}

\numberwithin{equation}{section}

\newtheorem{theorem}{Theorem}

\newtheorem{proposition}{Proposition}

\newtheorem{remark}{Remark}

\numberwithin{theorem}{section}
\numberwithin{corollary}{section}
\numberwithin{lemma}{section}
\numberwithin{definition}{section}
\numberwithin{proposition}{section}
\numberwithin{remark}{section}

\newcommand{\dint}{\ds\int}
\newcommand{\ds}{\displaystyle}

\newcommand{\rn}{\mathbb R^n}
\newcommand{\R}{\mathbb R}

\def\Sn{\mathbb{S}^n}

\newcommand{\medint}{-\kern  -,375cm\int}

\def\det{\mbox{det}\,}
\def\de{\partial}

\begin{document}
\title[ ]{Characterization of ellipsoids through an overdetermined boundary value problem of Monge-Amp\`ere type}
\author[B. Brandolini, N. Gavitone, C. Nitsch, C. Trombetti]{B. Brandolini$^*$ - N.
Gavitone$^*$ - C. Nitsch$^*$ - C. Trombetti$^*$}
\thanks{$^*$ Dipartimento di Matematica e Applicazioni ``R. Caccioppoli'',
Universit\`{a} degli Studi di Napoli ``Federico II'', Complesso Monte S.
Angelo, via Cintia - 80126 Napoli, Italy; email: brandolini@unina.it; nunzia.gavitone@unina.it; c.nitsch@unina.it; cristina@unina.it}

\begin{abstract}
The study of the optimal constant in an Hessian-type Sobolev inequality leads to a fully nonlinear boundary value problem, overdetermined with non standard boundary conditions. We show that all the solutions have ellipsoidal symmetry. In the proof we use the maximum principle applied to a suitable auxiliary function in conjunction with an entropy estimate from affine curvature flow.\end{abstract}
\maketitle

\section{Introduction}

In this paper we study the following fully nonlinear overdetermined boundary value problem
 \begin{equation}
\label{prob}
  \left\{
    \begin{array}{ll}
      \det D^2u=1 &\text{in } \Omega\\\\
      u=0 &\text{on } \de \Omega\\\\
      H_{n-1}|Du|^{n+1}=c &\text{on } \de \Omega,
    \end{array}
  \right.
\end{equation}
where $\Omega\subset \rn$ is a smooth, bounded open set whose boundary has positive Gaussian curvature $H_{n-1}$, and $c$ is a given positive constant. If we denote by $\omega_n$ the volume of the unit ball in $\R^n$ and $\Omega$ is any ellipsoid of measure $\omega_n c^{n/2}$, then 
\begin{equation}\label{ell}
u(x)=\dfrac{|A(x-x_0)|^2-c}{2}
\end{equation}
 is the solution to \eqref{prob}, for some $x_0\in \rn$ and some $n \times n$ matrix $A$ with $\det A=1$. Obviously   $\Omega=\left\{x \in \rn:\> |A(x-x_0)|^2 \le c\right\}$.
Our main result reads as follows.
\begin{theorem}
\label{main}
Let  $\Omega\subset \rn$ be a bounded, convex, open set with $C^2$ boundary; a convex function $u \in C^2(\bar{\Omega})$ is a solution to problem \eqref{prob} if and only if $\Omega$ is an ellipsoid of measure $\omega_n c^{{n}/{2}}$ and $u$ has the specific form given in \eqref{ell}.
\end{theorem}

In 1971 in a celebrated paper (\cite{s71}) Serrin proved  that a smooth domain $\Omega$ is necessarily a ball if, for some constant $\gamma >0$, there exists a solution $u \in C^2(\bar \Omega)$ to the following problem
\begin{equation}\label{eq_serrin}
 \left\{
    \begin{array}{ll}
      \Delta u=1 &\text{in } \Omega\\\\
      u=0 &\text{on } \de \Omega\\\\
      \dfrac{\partial u}{\partial \nu}=\gamma &\text{on } \de \Omega,
    \end{array}
  \right.
\end{equation}
where $\nu$ is the unit outer normal to $\de \Omega$.
The main ingredients employed in the proof were a revisited Alexandrov moving plane method and a refinement of the maximum principle and Hopf's boundary point Lemma. All such techniques soon became primary tools in the study of symmetries in PDE's (see for instance \cite{f00} and the references therein) when, in the wake of this pioneering paper, the study of overdetermined boundary value problems burst out.

Right after Serrin's paper the very same result was also obtained by Weinberger \cite{w71} with a very short proof. To better understand the key steps of our proof in the following Sections, it is worth to briefly remind here the basic ideas behind Weinberger's one. First of all he showed that the auxiliary function $|Du|^2-\frac{2}{n}u$ (being subharmonic in $\Omega$) achieves its maximum $\gamma^2$ on the boundary of $\Omega$. Then he observed that, in view of the Poho\v zaev identity, one has $$\int_\Omega |Du|^2\,dx-\frac{2}{n}\int_\Omega u\,dx=\gamma^2|\Omega|$$
 ($|\Omega|$ denoting the measure of $\Omega$), and he deduced that $|Du|^2-\frac{2}{n}u$ is constant in $\Omega$. This fact immediately carries the radial symmetry of the solution to \eqref{eq_serrin}. 

Since these fundamental contributions, several alternative proofs and generalizations to linear and nonlinear operators followed  (see for instance \cite{w71,gl89,ch98,hp98,bh02,fgk06,bnst,bnst08,fk08,bnst09s,bk11,se11,f12}). Maximum principle is always hidden somewhere in the proof, however some of the developed techniques do not require its explicit usage (we refer the interested reader to \cite{bnst,bnst08,bnst09s}). 

Compared to most of the problems that can be found in literature, \eqref{prob} has some unusual peculiarities. Firstly the differential operator is fully nonlinear, with strongly coupled second order derivatives. Secondly the problem admits both radially and non radially symmetric solutions. Such two features can be found in literature for instance in \cite{gnn79,l81,pp96,hp98, ps01,fbk03,c04,ap06,cs09,b13,dpg_mn}, where they rarely occur simultaneously and, to our knowledge, not for all dimensions. 


The structure of our paper is the following. In Section 2 we introduce basic notation and preliminary results. Section 3 is the core of the paper and for the reader's convenience we split the proof in four claims. In the wake of Weinberger's paper we introduce an auxiliary function $\varphi(u,Du,D^2 u)$ for which a maximum principle holds (see Claim 1 and Claim 2 below). In view of a Poho\v zaev type identity for Monge-Amp\`ere equations we show that $\varphi$ is constant in $\Omega$ (see Claim 3 below). Surprisingly, this provides informations on the evolution of $\de \Omega$ by affine mean curvature flow. In particular, an equality sign is achieved in a fundamental entropy inequality (involving the affine surface area of $\Omega$) which have been proved in \cite{and06} and as a consequence $\Omega$ turns out to be an ellipsoid (see Claim 4 below).

The use of the affine mean curvature machinery is somehow the most original idea in our proof. We observe that, at least in the planar case, such an idea is not needed (see \cite{ap06}), and for completeness we sketch a different proof in Remark \ref{rem2}.


\smallskip

Now, before enter in the details of the proof of Theorem \ref{main}, we want to discuss the reasons which led us to consider the overdetermination in \eqref{prob}. 
They have to be found in connection with the study of Hessian Sobolev inequalities. To better understand such a link we give an insight on how classical isoperimetric estimates on the best constant in a Sobolev-Poincar\'e type inequality are related to classical overdetermined problems like \eqref{eq_serrin}.
To this aim assuming that $\Omega$ is a bounded, open subset of $\R^n$ then there exists the least positive constant $\mathcal{T}(\Omega)$ (called torsional rigidity of $\Omega$) such that
$$\left(\int_\Omega u\,dx\right)^2\le \mathcal{T}(\Omega)\int_{\Omega}|Du|^2\,dx,$$ 
for all $u\in H_0^1(\Omega).$
The variational characterization of the torsional rigidity is 
\begin{equation}\label{eq_vart}
\frac{1}{\mathcal{T}(\Omega)}=\min_{v \in H_0^1(\Omega)} \frac{\dint_\Omega |Dv|^2 \,dx}{\left(\dint_\Omega v\,dx\right)^{2}},
\end{equation}
and any function achieving the minimum on the right-hand side of \eqref{eq_vart} is proportional to the unique solution to the following Poisson problem
\begin{equation*}\label{eq_poisson}
\left\{\begin{array}{ll}
\Delta u=1 & \mbox{in}\> \Omega\\
\\
u=0 &\mbox{on}\> \partial\Omega.
\end{array}
\right.
\end{equation*}
Now, under suitable smoothness assumptions on $\Omega$, it is possible to differentiate the torsional rigidity $\mathcal{T}(\Omega)$ with respect to any smooth domain variation, leading to Hadamard formula (see for instance \cite{had,henrot-pierre}). Under the additional constraint of keeping the measure of $\Omega$ fixed, we can call stationary domains those smooth open sets on which the derivatives of $\mathcal{T}(\cdot)$ along any smooth domain deformation vanish. It turns out that a domain is stationary if it admits a solution to problem \eqref{eq_serrin} and, according to such a notion, Serrin has proved that balls are the unique stationary domains.  
This is somehow in agreement with classical results (see for instance \cite{polya,ta}) which established that, among sets of given measure, $\mathcal{T}(\Omega)$ is maximal on balls (observe that under the same prescription $\mathcal{T}(\Omega)$ is not bounded away from zero). This property was first noticed by a famous mechanician of the 19th century and named after him {\it Saint-Venant's Principle}. Serrin's paper strengthened such a principle by proving that no other stationary domain exists.

We turn now our attention to higher order (Hessian-type)  Sobolev inequalities.
In the paper \cite{cw01} (see also \cite{t90,w94,choge,truwan,ga1}), among other things, the authors proved that,  whenever $\Omega$ is a smooth, convex set, there exists the least positive constant $\mathcal{S}(\Omega)$ such that the following Hessian Sobolev inequality holds
\begin{equation*}\label{sobolev}
\left(\int_\Omega (-u) dx \right)^{n+1} \le \mathcal{S}(\Omega) \int_\Omega(-u)\det D^2 u\,dx
\end{equation*}
for all $u\in \Phi_0(\Omega)\equiv\left\{u\in C^2(\Omega)\cap C^0(\bar \Omega):\> u \>\mbox{convex in}\> \Omega, \> u=0 \>\mbox{on} \>\partial\Omega\right\}$. The value of $\mathcal{S}(\Omega)$ 
can be characterized by
\begin{equation}\label{best-sobolev}
\frac{1}{\mathcal{S}(\Omega)}=\min_{w \in \Phi_0(\Omega)} \frac{\dint_\Omega (-w)\det D^2 w \,dx}{\left(\dint_\Omega (-w)\,dx\right)^{n+1}}
\end{equation}
and any function in $\Phi_0(\Omega)$ achieving the minimum on the right-hand side of \eqref{best-sobolev} is proportional to the unique solution to the following Monge-Amp\`ere boundary value problem
\begin{equation}\label{b}
\left\{\begin{array}{ll}
\det D^2 u=1 & \mbox{in}\> \Omega\\
\\
u=0 &\mbox{on}\> \partial\Omega.
\end{array}
\right.
\end{equation}
We can call $S(\Omega)$ the Monge-Amp\`ere torsional rigidity in analogy with the definition given above.
Now, as for the classical torsional rigidity, we want to identify the stationary domains. To this aim let us consider a smooth, strictly convex open set $\Omega$ and a family of maps $\Psi(t)$ satisfying
$$
\Psi: t \in [0,T[ \to W^{1,\infty}(\rn,\rn)\> \mbox{differentiable at }\>0\> \>\mbox{with } \>\Psi(0)=I, \> \Psi'(0)=V
$$
where $I$ is the identity and $V$ is a vector field. Let us denote $\Omega_t=\Psi(t)(\Omega)$ and 
$$
s(\Omega_t)=\mathcal{S}(\Omega_t)^{1/n}=\int_{\Omega_t} (-u(x,t))\,dx 
$$
where $u(x,t)$ solves problem \eqref{b} with $\Omega$ replaced by $\Omega_t$. By Hadamard formula  we get 
$$
s'(\Omega_t)\Big|_{t=0}=\frac{1}{n}\int_{\partial\Omega} H_{n-1}|Du|^{n+1} V\cdot \nu.
$$
Thus $\Omega$ is a stationary point for $\mathcal{S}(\Omega)$ under the volume constraint if 
$$
\int_{\partial\Omega} H_{n-1}|Du|^{n+1} V\cdot \nu=0, \quad \mbox{for any} \>V \>\mbox{such that} \int_{\partial\Omega}V \cdot \nu=0.
$$
Therefore a stationary domain $\Omega$ carries the additional condition $$H_{n-1}|Du|^{n+1}=const$$ and for such a set a solution to problem \eqref{prob} exists. 
In this framework we can read Theorem \ref{main} as the proof that no stationary domain for $\mathcal{S}(\cdot)$ exists other than ellipsoids. Our result is in agreement with previous papers \cite{bnt09, bnt11} where it has been proved that in the class of smooth, strictly convex, open sets of given measure, $\mathcal{S}(\cdot)$ is minimal on all ellipsoids. Balls are not the only domains since the Monge-Amp\`ere operator is invariant under measure preserving affine transformations and therefore it is unable to ``distinguish'' a ball from an ellipsoid. 
The analogy between $\mathcal{S}(\cdot)$ and $\mathcal{T}(\cdot)$ is not as tight as it might seem since, contrary to the behavior  of $\mathcal{T}(\cdot)$, the constant $\mathcal{S}(\cdot)$ happens to be minimal on balls and not maximal. However, what it really makes a difference, is that once  $\mathcal{S}(\cdot)$ is continuously extended to the whole class of convex sets,  trivial arguments involving maximum principle ensure that such a constant is also bounded from above in terms of the measure of $\Omega$ alone. Compactness results in the class of convex sets (Blaschke-Santal\`o theorem) guarantee that the maximum is achieved.  
The determination of maximizers is a puzzling nontrivial open problem. 
As a corollary to our result we deduce that 
the maximum of $\mathcal{S}(\cdot)$ is achieved on sets which are convex but do not belong to the class of $C^2$ strictly convex sets.

\section{Notation and preliminaries}

\subsection{Symmetric functions and Hessian operators}
We denote by $A=(a_{ij})$ a matrix in the space $\Sn$ of the real symmetric $n\times n$ matrices, and by
$\lambda_1,...,\lambda_n$ its eigenvalues. For $k\in
\left\{1,...,n\right\}$, the $k$-th elementary symmetric function
of $A$ is
$$S_k(A)=S_k(\lambda_1,...,\lambda_n)=\sum_{1\leq i_1
<\cdots<i_k\leq n } \lambda_{i_1}\cdots \lambda_{i_k}.$$ 
Note that $S_k(A)$ is just the sum of all $k \times k$ principal minors of $A$.

The
operator $S_k^{1/k}$, for $k=1,...,n$, is homogeneous of degree
$1$ and it is concave, if restricted to 
$$\Gamma_k=\{A\in\Sn\,:\,S_i(A)\geq 0\text{ for }i=1,\dots,k\}\,.$$
Denoting by
$$
S_k^{ij}(A)
= \frac{\partial }{\partial a_{ij}}S_k(A),$$
Euler identity for homogeneous functions gives
$$
S_k(A) = \frac{1}{k} \sum_{i,j}S_k^{ij}(A) a_{ij}.
$$ 

We will use the following notations:
\begin{itemize} 
\item $S_k(i)$ means the $k$-th elementary symmetric function of $\lambda_1,...,\lambda_n$ excluding $\lambda_i$; 

\item $S_k(i,j)$ means the $k$-th elementary symmetric function of $\lambda_1,...,\lambda_n$ excluding $\lambda_i$ and $\lambda_j$; 

\item $S_k(i,j,r)$  means the $k$-th elementary symmetric function of $\lambda_1,...,\lambda_n$ excluding $\lambda_i$, $\lambda_j$ and $\lambda_r$;

\item finally, when $k=n$,
\begin{equation}\label{derdet}
S_n^{ij}(A)=\frac{\de (\det A)}{\de a_{ij}},\qquad S_n^{ij,rs}(A)=\frac{\de^2 (\det A)}{\de a_{ij}\de a_{rs}}, \qquad  
S_n^{ij,rs,\alpha \beta}(A)=\frac{\de^3 (\det A)}{\de a_{ij}\de a_{rs} \de a_{\alpha \beta}}. 
\end{equation}
\end{itemize}

If $A$ has diagonal form, \eqref{derdet} becomes (see for example \cite{gm03})
\begin{gather}
\label{derdiag1}
\begin{split} 
S_n^{ij}(A)&=\begin{cases}
\begin{array}{ll}
  S_{n-1}(i) & \text{ if } i=j \\
   0 &  \text{ if } i\neq j
\end{array}
\end{cases} \\
 S_n^{ij,rs}(A)&=\begin{cases}
\begin{array}{ll}
  S_{n-2}(i,j) & \text{ if } i=j, r=s, i \neq r \\
 -S_{n-2}(i,j) &  \text{ if } i\neq j, r=j, s=i\\
 0 & \text{ otherwise}
\end{array}
\end{cases} 
\\
S_n^{ij,rs,\alpha \beta}(A)&=\begin{cases}
\begin{array}{ll}
  S_{n-3}(i,r,\alpha) & \text{ if } i=j, r=s, \alpha=\beta, r \neq i, \alpha \neq i,r \\
 -S_{n-3}(i,r, \alpha) &  \text{ if } i= j, r\neq s, \alpha=s, \beta=r,   r \neq i, \alpha \neq i,r\\
-S_{n-3}(i,r, \alpha) & \text{ if } i\neq j, r=s, \alpha=\beta, r \neq i, \alpha \neq i,r\\
S_{n-3}(i,r,\alpha) & \text{ if } i\neq j, r=j, s=\alpha, \beta=i,   r \neq i, \alpha \neq i,r\\
 0 & \text{ otherwise}. 
\end{array}
\end{cases}
\end{split}
\end{gather}

Now let $\Omega$ be an open subset of $\rn$ and let $u\in C^2(\Omega)$. The $k$-Hessian operator $S_k\left(D^2u\right)$
is defined as the $k$-th elementary symmetric function of $D^2u$.
Notice that $$S_1(D^2u)=\Delta u\quad \mbox{and}\quad
S_n(D^2u)=\mbox{det}D^2u.$$
For $k>1$, the $k$-Hessian operators are fully nonlinear and, in general, not elliptic,
unless restricted to the class of $k$-convex functions
$$\Phi_k^2(\Omega)=\left\{u\in C^2(\Omega)\,:\,S_i(D^2u)\geq 0 \text{
in }\Omega, i=1,2,...,k \right\}.$$ Notice that
$\Phi_n^2(\Omega)$ coincides with the class of $C^2(\Omega)$
convex functions.

\noindent A direct computation yields that $(S_k^{1j}(D^2u),\dots,S_k^{nj}(D^2u))$ is divergence free, i.e.
\begin{equation}
\label{div0}
\sum_i\frac{\partial}{\partial x_i}S_k^{ij}=0, \qquad j=1,...,n;
\end{equation}
hence $S_k(D^2 u)$ can be written in divergence form
\begin{equation}\label{divk}
S_k(D^2 u) = \frac{1}{k} S_k^{ij}(D^2 u) u_{ij} =
 \frac{1}{k} (S_k^{ij}(D^2 u) u_j)_i,
\end{equation}
where subscripts stand for partial differentiations.

If $t$ is a regular value of $u$ and $H_{n-1}$ stands for the Gaussian curvature of the level set $\de\{u \le t\}$ at the point $x$, the following pointwise identity holds (see \cite{re74}, \cite{tr97})
\begin{equation}\label{H-S}
H_{n-1}=\frac{S_n^{ij}(D^2 u)u_iu_j}{|Du|^{n+1}}.
\end{equation}

Finally we recall the following Poho\v zaev identity (see \cite{tso90}, \cite{bnst}) 
\begin{proposition}
Let  $\Omega \subset \R^n$ be a bounded, convex, open set with $C^2$ boundary and let $f\in C^1(\R)$ be a nonnegative function. If $u \in C^2(\bar{\Omega})$ is a convex solution to the problem  
\begin{equation*}
  \left\{
    \begin{array}{ll}
      \mathrm{det}\, D^2 u=f(u) &\text{in } \Omega\\\\
      u=0 &\text{on } \de \Omega,
    \end{array}
  \right.
\end{equation*}
then the following identity holds
\begin{equation}\label{poz}
-\frac{1}{n+1} \int_{\Omega}S_n^{ij}(D^2 u)u_i u_j \,dx+\frac{1}{n+1} \int_{\de \Omega} \langle x, \nu\rangle H_{n-1} |Du|^{n+1}=n \int_{\Omega}F(u)\,dx,
\end{equation}
where $S_n^{ij}(D^2 u)$ are defined in \eqref{derdet}, $\nu$ is the outer unit normal to $\de \Omega$,  and $F(u)=\int_u^0 f(s)ds$.   
\end{proposition}

\subsection{The affine curvature flow}
An affine isoperimetric inequality, known as Petty inequality, states that (see, for example, \cite{petty, schneider})
\begin{equation}\label{petty}
   \int_{\partial \Omega}H_{n-1}^{ \frac{1}{n+1}}\le n  \omega_n^{\frac{2}{n+1}} |\Omega|^{\frac {n-1}{n+1}},
\end{equation}
 equality holding if and only if $\Omega$ is an ellipsoid. The integral on the left hand side of \eqref{petty} is known as affine surface area and $H_{n-1}^{ \frac{1}{n+1}}$ is known as affine curvature.
 
Now we recall a result proved in \cite{and06} concerning affine curvature flow. 
\begin{theorem}[Andrews '96]\label{lemma_Andrews}
Let $\varphi_0:{\mathbb S}^{n-1}\to \R^n$ be a smooth ($C^\infty$) strictly convex embedding of the unit sphere in  $\R^n$; then there exists a unique $t_E>0$ and a unique 
$$\varphi(z,t)\in C^\infty({\mathbb S}^{n-1} \times  [0,t_E[;\R^n)$$ such that, for all 
 $0\le t<t_E, $ $\varphi(\cdot, t):
{\mathbb S}^{n-1}\to \Gamma_t \subset \R^n$ is a smooth, closed surface, uniformly convex (i.e. with  strictly positive Gaussian curvature) for $t>0$,   
and for all $(z,t) \in {\mathbb S}^{n-1} \times [0,t_E[$ $\varphi$ is a
solution to the following partial differential equation
\begin{equation}\label{evol}
\frac{\partial \varphi}{\partial t}(z,t)=-(H_{n-1}[\Gamma_t](\varphi(z,t)))^{\frac{1}{n+1}} \nu_{\Gamma_t}(\varphi(z,t)),
\end{equation}
where $H_{n-1}[\Gamma_t](x)$ and $\nu_{\Gamma_t}(x)$ are, respectively, the Gaussian curvature and the outer normal of $\Gamma_t$ at the point $x\in\Gamma_t$ and $\varphi(z, 0) = \varphi_0(z)$.

Moreover 
\begin{enumerate}[i)]
\item $\Gamma_t$ converges to a point as $t\nearrow t_E,$
\item after rescaling about the final point to make the enclosed volume constant, $\Gamma_t$ converges in $C^\infty$ to an ellipsoid,
\item the following estimate holds:
\begin{equation}
\label{derandr}
\frac{\de}{\de t}\left(V_t^{-\frac{n-1}{n+1}}\int_{\Gamma_t}H_{n-1}[\Gamma_t]^{\frac{1}{n+1}}\right)\ge 0,
\end{equation}
and the inequality is strict unless $\varphi({\mathbb S}^{n-1},t)$ is an ellipsoid  for any $0 \le t < t_E $. 
Here $V_t$ is the enclosed volume of the hypersurface $\varphi({\mathbb S}^{n-1},t)$.
\end{enumerate} 
\end{theorem}

\noindent
The affine curvature flow of a convex surface is a flow where each point of the surface moves in the direction of the inner normal with velocity equal to the affine curvature of the surface itself. The previous theorem states that, for any initial smooth, convex, closed surface,  it is possible to find a unique one parameter family of solutions to the affine curvature flow. Such a family is smooth and shrinks to a point by approaching an ellipsoidal shape.

\section{Proof of the theorem \ref{main}}

From now on $u\in C^2(\bar \Omega)\cap C^\infty(\Omega)$ will be a convex solution to problem \eqref{prob}. For the reader's convenience we will denote
\begin{equation*}\label{derdet2}
S^{ij}=S_n^{ij}(D^2 u),\qquad S^{ij,rs}=S_n^{ij,rs}(D^2 u), \qquad  
S^{ij,rs,\alpha\beta}=S_n^{ij,rs,\alpha \beta}(D^2 u). 
\end{equation*}
Differentiating the equation in \eqref{prob} we immediately get that, for $k=1,...,n$, 
\begin{equation}\label{der_det}
\frac{\partial (\det D^2u)}{\partial x_k}= \sum_{i,j=1}^n S^{ij}u_{ijk}=0 \qquad \mathrm{in}\> \Omega,
\end{equation}
and, for $k,l=1,...,n$,
\begin{equation}\label{der_sec_det}
\frac{\partial^2 (\det D^2u)}{\partial x_k \de x_l}= \sum_{\substack{i,j\\r,s}} S^{ij,rs}u_{ijk}u_{rsl}+\sum_{i,j} S^{ij}u_{ijkl}=0  \qquad \mathrm{in}\> \Omega.
\end{equation}
For any $v \in C^2(\Omega)$,  let us consider the linear operator
\begin{equation*}
\label{oplin}
Lv=\sum_{i,j} S^{ij} v_{ij},
\end{equation*}
where $(S^{ij})$ is the cofactor matrix of $D^2 u$. Then, for any $k=1,...,n$, by \eqref{der_det} we immediately get 
\begin{equation}\label{eq_partial}
Lu_k=0 \qquad \mbox{in}\> \Omega.
\end{equation}
Let us introduce the following auxiliary function
\begin{equation}
\label{fi}
\varphi=H_{n-1}|Du|^{n+1}-2u.
\end{equation}

\medskip
\noindent {\bf Claim 1: the following identity holds true
\begin{equation} 
\label{L_fine}
L\varphi= \sum_{\substack{k,l\\i,j}}\sum_{\substack{r,s\\ \alpha,\beta}} S^{kl} S^{ij, rs, \alpha \beta}u_{rsk} u_{\alpha \beta l}u_i u_j-\sum_{\substack{k,l\\i,j}}\sum_{\substack{r,s\\ \alpha,\beta}} S^{ij, kl} S^{rs,\alpha \beta}u_{rsk}u_{\alpha \beta l} \,u_i u_j.
\end{equation}}
Being $\det D^2 u=1$, the matrix $(S^{ij})$ is the inverse of $D^2 u$, that is $S^{ij} u_{ik}=\delta_{j k}$. Hence by \eqref{H-S} we get
\begin{eqnarray*}
\frac{\de}{\de x_k}\left( H_{n-1} |Du|^{n+1}\right)&=& \frac{\de}{\de x_k}\left( \sum_{i,j} S^{ij}u_i u_j\right)= \sum_{\substack{i,j\\r,s}} S^{ij, rs}u_{rsk}u_i u_j+2 \sum_{i,j} S^{ij}u_{ik} u_j
\\
&=& \sum_{\substack{i,j\\r,s}} S^{ij, rs}u_{rsk}u_i u_j+2 u_{k}, 
\end{eqnarray*}
and then
\begin{equation*}
\frac{\de \varphi}{\de x_k} =\sum_{\substack{i,j\\r,s}} S^{ij, rs}u_{rsk}u_i u_j.
\end{equation*}
Let us compute the second derivatives of $\varphi$
\begin{equation*}
\frac{\de^2 \varphi}{\de x_k \de x_l} =\sum_{\substack{i,j\\r,s}} \frac{\de}{\de x_l}\left(S^{ij, rs}\right)u_{rsk}u_i u_j+\sum_{\substack{i,j\\r,s}} S^{ij, rs}u_{rskl}u_i u_j+ 2\sum_{\substack{i,j\\r,s}} S^{ij, rs}u_{rsk}u_{i l} u_j.
\end{equation*}
Hence
\begin{gather}
\label{espr}
\begin{split}
L \varphi &=\sum_{k,l}S^{kl}\varphi_{kl}
\\
&=\sum_{k,l}\sum_{\substack{i,j\\r,s}} S^{kl}\frac{\de}{\de x_l}\left(S^{ij, rs}\right)u_{rsk}u_i u_j+\sum_{k, l}\sum_{\substack{i,j\\r,s}} S^{kl} S^{ij, rs}u_{rskl}u_i u_j+ 2\sum_{k, l}\sum_{\substack{i,j\\r,s}} S^{kl} S^{ij, rs}u_{rsk}u_{i l} u_j \\
&=\sum_{\substack{i,j\\k,l}}\sum_{\substack{r,s\\\alpha, \beta}} S^{kl} S^{ij, rs, \alpha \beta}u_{rsk} u_{\alpha \beta l}u_i u_j+\sum_{k, l}\sum_{\substack{i,j\\r,s}} S^{kl} S^{ij, rs}u_{rskl}u_i u_j + 2\sum_{\substack{i,j\\r,s}}  S^{ij, rs}u_{rsi} u_j.
\end{split}
\end{gather}
By \eqref{div0} we get that the last term in \eqref{espr} vanishes and then 
\begin{equation*}
L \varphi =\sum_{\substack{i,j\\k,l}}\sum_{\substack{r,s\\\alpha, \beta}} S^{kl} S^{ij, rs, \alpha \beta}u_{rsk} u_{\alpha \beta l}u_i u_j+\sum_{k, l}\sum_{\substack{i,j\\r,s}} S^{kl} S^{ij, rs}u_{rskl}u_i u_j.
\end{equation*}
Finally substituting \eqref{der_sec_det} in the above equality we get \eqref{L_fine}.

\medskip
\noindent {\bf Claim 2:  $L\varphi \ge 0 \text{ in }\Omega.$}

\noindent We distinguish two cases. Suppose first that $n=2$. We observe that \eqref{eq_partial} gives

\begin{equation}\label{alphabeta}
\det D^2(\alpha u_1+\beta u_2)=\det\left(\alpha D^2 u_1+\beta D^2u_2\right)\le 0 \qquad \forall \alpha, \beta\in \R.
\end{equation}
Therefore we immediately have
\begin{eqnarray}\label{eq2}
L\varphi&=&-\sum_{\substack{k,l\\i,j}}\sum_{\substack{r,s\\ \alpha,\beta}} S^{ij, kl} S^{rs,\alpha \beta}u_{rsk}u_{\alpha \beta l} \,u_i u_j\\
&=&2u_1^2(-\det D^2 u_2)+2u_2^2(-\det D^2 u_1)-2u_1u_2(u_{111}u_{222}-u_{112}u_{122})\notag
\\
&=& -2\det\left(u_2D^2 u_1-u_1D^2u_2\right)\ge 0.\notag
\end{eqnarray} 

Suppose now that $n>2$. Let $x \in \Omega$; by performing a rotation of the coordinates we may suppose that $D^2 u(x)$ is diagonal. Since $\det D^2 u=1$, \eqref{derdiag1} can be rewritten as follows
\begin{gather}
\label{derdiagl}
\begin{split} 
S^{ij}&=\begin{cases}
\begin{array}{ll}
  \frac{1}{\lambda_i} & \text{ if } i=j \\
   0 &  \text{ if } i\neq j
\end{array}
\end{cases} \\
 S^{ij,rs}&=\begin{cases}
\begin{array}{ll}
  \frac{1}{\lambda_i\lambda_j} & \text{ if } i=j, r=s, i \neq r \\
 -\frac{1}{\lambda_i\lambda_j} &  \text{ if } i\neq j, r=j, s=i\\
 0 & \text{ otherwise} 
\end{array}
\end{cases} 
\\
S^{ij,rs,\alpha \beta}&=\begin{cases}
\begin{array}{ll}
  \frac{1}{\lambda_i\lambda_r\lambda_\alpha} & \text{ if } i=j, r=s, \alpha=\beta, r \neq i, \alpha \neq i,r \\
 - \frac{1}{\lambda_i\lambda_r\lambda_\alpha}  &  \text{ if } i= j, r\neq s, \alpha=s, \beta=r,   r \neq i, \alpha \neq i,r\\
- \frac{1}{\lambda_i\lambda_r\lambda_\alpha}  & \text{ if } i\neq j, r=s, \alpha=\beta, r \neq i, \alpha \neq i,r\\
 \frac{1}{\lambda_i\lambda_r\lambda_\alpha}  & \text{ if } i\neq j, r=j, s=\alpha, \beta=i,   r \neq i, \alpha \neq i,r\\
 0 & \text{ otherwise.} 
\end{array}
\end{cases}
\end{split}
\end{gather}
First let us consider the second term in the right-hand side of \eqref{L_fine}. By using \eqref{derdiagl} we have 
\begin{gather}
\begin{split}\label{pezzo2}
&\sum_{\substack{k, l\\ i,j}}\sum_{\substack{r,s\\ \alpha,\beta}} S^{ij, kl} S^{rs,\alpha \beta}u_{rsk}u_{\alpha \beta l} \,u_i u_j
\\
&= \left(\sum_{i=j}+\sum_{i\neq j} \right)\sum_{k,l}\sum_{\substack{r,s\\ \alpha,\beta}} S^{ij, kl} S^{rs,\alpha \beta}u_{rsk}u_{\alpha \beta l} \,u_i u_j 
\\
&=\sum_i\sum_{k \neq i}u_i^2S^{ii,kk} \sum_{\substack{r,s\\ \alpha,\beta}} S^{rs,\alpha \beta}u_{rsk}u_{\alpha \beta k}+ \sum_i \sum_{j\neq i}u_i u_j S^{ij,ji}\sum_{\substack{r,s\\ \alpha,\beta}} S^{rs,\alpha \beta}u_{rsj}u_{\alpha \beta i}
\\
&=\sum_i\sum_{k \neq i}u_i^2S^{ii,kk} \sum_r\sum_{\alpha\neq r}\left[ S^{rr,\alpha \alpha}u_{rrk}u_{\alpha \alpha k}+
S^{r\alpha,\alpha r}u^2_{r\alpha k}\right] 
\\
&\quad+ \sum_i \sum_{j\neq i} u_i u_j S^{ij,ji}\sum_r\sum_{\alpha\neq r}\left[ S^{rr,\alpha \alpha}u_{rrj}u_{\alpha \alpha i}+S^{r\alpha,\alpha r}u_{r\alpha j}u_{r\alpha i}\right]
\\
&=\sum_i\sum_{k \neq i}\frac{u_i^2}{\lambda_i \lambda_k}\sum_r\sum_{\alpha \neq r} \frac{1}{\lambda_r \lambda_\alpha}\left(u_{rrk}u_{\alpha \alpha k}- u^2_{r\alpha k}\right)- \sum_i \sum_{j\neq i}\frac{u_i u_j}{\lambda_i \lambda_j}\sum_r \sum_{\alpha \neq r}  \frac{1}{\lambda_r \lambda_\alpha}\left(u_{rrj}u_{\alpha \alpha i}-u_{r\alpha j}u_{r\alpha i}\right) 
\\
&=\sum_i \frac{1}{\lambda_i}\sum_{j\neq i} \frac{1}{\lambda_j} \sum_r\frac{1}{\lambda_r}\sum_{\alpha \neq r}\frac{1}{\lambda_\alpha}\left[u_i^2 \left(u_{rrj}u_{\alpha \alpha j}- u^2_{r\alpha j}\right) -u_i u_j\left(u_{rrj}u_{\alpha \alpha i}-u_{r\alpha j}u_{r\alpha i}\right)\right]. 
\end{split}
\end{gather}
Analogously, by \eqref{derdiagl} the first term in the right-hand side of \eqref{L_fine} becomes
\begin{gather}
\begin{split}
\label{pezzo1}
&\sum_{\substack{k,l\\i,j}}\sum_{\substack{r,s\\ \alpha,\beta}} S^{kl} S^{ij, rs, \alpha \beta}u_{rsk} u_{\alpha \beta l}u_i u_j
\\
&=\sum_k\sum_{\substack{i,j,r,s\\ \alpha,\beta}} S^{kk} S^{ij, rs, \alpha \beta}u_{rsk} u_{\alpha \beta k}\, u_i u_j 
\\
&= \left(\sum_{i=j}+\sum_{i\neq j} \right) \sum_k\sum_{\substack{r,s\\ \alpha,\beta}}  S^{kk} S^{ij, rs, \alpha \beta}u_{rsk} u_{\alpha \beta k}\, u_i u_j
\\
 &=\sum_{k} \sum_i\sum_{\substack{r,s\\ \alpha,\beta}}  S^{kk} S^{ii, rs, \alpha \beta}u_{rsk} u_{\alpha \beta k}\, u^2_i +\sum_{k} \sum_{i \neq j}\sum_{\substack{r,s\\ \alpha,\beta}} S^{kk} S^{ij, rs, \alpha \beta}u_{rsk} u_{\alpha \beta k}\, u_i u_j
  \\
   &=\sum_{k} \sum_iS^{kk}u_i^2\sum_{r \neq i}\sum_{\substack{\alpha\neq r,i}}\left(  S^{ii, rr, \alpha \alpha} u_{rrk} u_{\alpha \alpha k}+S^{ii,r\alpha,\alpha r} u^2_{r\alpha k}\right) 
   \\
&\quad+2 \sum_{k}\sum_i \sum_{j \neq i}S^{kk} u_i u_j\sum_{\substack{\alpha\neq i,j}}  \left(S^{ij, ji, \alpha \alpha}u_{ijk} u_{\alpha \alpha k}+S^{ij,j\alpha,\alpha i}u_{j \alpha k}u_{\alpha ik}  \right)
\\
&=\sum_{k} \frac{1}{\lambda_k}\sum_i \frac{1}{\lambda_i} \sum_{r \neq i}\frac{1}{\lambda_r}\sum_{\substack{\alpha\neq r,i}}\frac{1}{\lambda_\alpha}\left[u_i^2\left( u_{rrk} u_{\alpha \alpha k}- u^2_{r\alpha k}\right)+2u_i u_r\left( u_{r\alpha k}u_{\alpha i k}-u_{irk} u_{\alpha \alpha k}\right)\right].
\end{split}
\end{gather}
Joining \eqref{pezzo2} and \eqref{pezzo1} we get
\begin{equation*}
\begin{split}
L \varphi=&\sum_{k} \frac{1}{\lambda_k}\sum_i \frac{1}{\lambda_i} \sum_{r \neq i}\frac{1}{\lambda_r}\sum_{\substack{\alpha\neq r,i}}\frac{1}{\lambda_\alpha}\left[u_i^2\left( u_{rrk} u_{\alpha \alpha k}- u^2_{r\alpha k}\right)+2u_i u_r\left( u_{r\alpha k}u_{\alpha i k}-u_{irk} u_{\alpha \alpha k}\right)\right]+\\
&- \sum_k \frac{1}{\lambda_k}\sum_i \frac{1}{\lambda_i}\sum_{r\neq i} \frac{1}{\lambda_r} \sum_{\alpha \neq k}\frac{1}{\lambda_\alpha}\left[u_i^2 \left(u_{kkr}u_{\alpha \alpha r}- u^2_{k\alpha r}\right) -u_i u_r\left(u_{kkr}u_{\alpha \alpha i}-u_{k\alpha r}u_{k\alpha i}\right)\right].
\end{split}
\end{equation*}
We can rearrange the terms appearing in the above formula and we obtain 
\begin{gather*}
\begin{split}
L \varphi=&\sum_{k} \frac{1}{\lambda_k}\sum_i \frac{1}{\lambda_i} \sum_{r \neq i}\frac{1}{\lambda_r}\left\{\sum_{\substack{\alpha\neq r,i}}\frac{u_i^2}{\lambda_\alpha}\left( u_{rrk} u_{\alpha \alpha k}- u^2_{r\alpha k}\right) -\sum_{\alpha \neq k}\frac{u_i^2}{\lambda_\alpha}\left(u_{kkr}u_{\alpha \alpha r}- u^2_{k\alpha r}\right) \right\}\\
&+\sum_{k} \frac{1}{\lambda_k}\sum_i \frac{1}{\lambda_i} \sum_{r \neq i}\frac{1}{\lambda_r}\left\{\sum_{\substack{\alpha\neq r,i}}\frac{2u_i u_r}{\lambda_{\alpha}}\left( u_{r\alpha k}u_{\alpha i k}-u_{irk} u_{\alpha \alpha k}\right)+\sum_{\alpha \neq k}\frac{u_i u_r}{\lambda_\alpha}\left(u_{kkr}u_{\alpha \alpha i}-u_{k\alpha r}u_{k\alpha i}\right)\right\}\\
=& \sum_{k} \frac{1}{\lambda_k}\sum_i \frac{1}{\lambda_i} \sum_{r \neq i}\frac{1}{\lambda_r}\left(u_i^2\mathcal{I}_1+u_iu_r\mathcal{I}_2\right),
\end{split}
\end{gather*}
where 
\begin{equation*}
\mathcal{I}_1=\sum_{\substack{\alpha\neq r,i}}\frac{\left( u_{rrk} u_{\alpha \alpha k}- u^2_{r\alpha k}\right)}{\lambda_\alpha} -\sum_{\alpha \neq k}\frac{\left(u_{kkr}u_{\alpha \alpha r}- u^2_{k\alpha r}\right)}{\lambda_\alpha}
\end{equation*}
and
\begin{equation*}
\mathcal{I}_2=2\sum_{\substack{\alpha\neq r,i}}\frac{\left( u_{r\alpha k}u_{\alpha i k}-u_{irk} u_{\alpha \alpha k}\right)}{\lambda_{\alpha}}+\sum_{\alpha \neq k}\frac{\left(u_{kkr}u_{\alpha \alpha i}-u_{k\alpha r}u_{k\alpha i}\right)}{\lambda_\alpha}.
\end{equation*}
Let us first consider $\mathcal{I}_1$:
\begin{eqnarray*}
&\mathcal{I}_1=\displaystyle\sum_{\substack{\alpha \neq i,r,k}}\frac{\left( u_{rrk} u_{\alpha \alpha k}- u_{kkr}u_{\alpha \alpha r}\right)}{\lambda_\alpha}+ \frac{\left(u_{rrk}u_{kkk}-u^2_{rkk}\right)}{\lambda_k}-\frac{\left(u_{kkr}u_{rrr}-u^2_{krr}\right)}{\lambda_r}- \frac{\left(u_{kkr}u_{iir}-u^2_{kir}\right)}{\lambda_i}\qquad\qquad\qquad
\\
&=\left\{u_{rrk}\left( \displaystyle\sum_{\substack{\alpha \neq i,r,k}}\frac{ u_{\alpha \alpha k}}{\lambda_\alpha}\right)- u_{kkr}\left(\displaystyle\sum_{\substack{\alpha \neq i,r,k}}\frac{u_{\alpha \alpha r}}{\lambda_\alpha}\right)
+\displaystyle \frac{\left(u_{rrk}u_{kkk}-u^2_{rkk}\right)}{\lambda_k}\-\frac{\left(u_{kkr}u_{rrr}-u^2_{krr}\right)}{\lambda_r}- \frac{\left(u_{kkr}u_{iir}-u^2_{kir}\right)}{\lambda_i} 
\right\}.
\end{eqnarray*}
Using \eqref{der_det} we have
\begin{equation}
\label{rel}
\sum_{\substack{\alpha\neq i,r,k}}\frac{u_{\alpha \alpha \beta}}{\lambda_{\alpha}}=-\left(\frac{u_{ii\beta}}{\lambda_i}+\frac{u_{rr\beta}}{\lambda_r}+\frac{u_{kk\beta}}{\lambda_k} \right), \qquad \text{for }\beta=1,\ldots,n,
\end{equation} 
and hence
\begin{equation}
\label{fine1}
\mathcal{I}_1=-\frac{\left( u_{rrk} u_{iik}-u^2_{kir}\right)}{\lambda_i}. 
\end{equation}
Reasoning in an analogous way as before we get
\begin{equation}
\label{fine2}
\mathcal{I}_2=\sum_{\alpha}\frac{u_{r \alpha k} u_{\alpha i k}}{\lambda_{\alpha}} .
\end{equation}
Using \eqref{fine1} and \eqref{fine2} we can write
\begin{gather*}
\begin{split}
L \varphi&=\sum_k \sum_i \sum_{r \neq i} \sum_{\alpha}\frac{u_i u_r}{\lambda_{\alpha}\lambda_i \lambda_k \lambda_r} u_{r \alpha k} u_{\alpha i k}-\sum_k \sum_i \sum_{r \neq i}\frac{u_i^2}{\lambda_k \lambda_r\lambda^2_i}\left( u_{rrk} u_{iik}-u^2_{kir}\right)\\
&=\sum_k \sum_i \sum_{r \neq i} \sum_{\alpha}\frac{u_i u_r}{\lambda_{\alpha}\lambda_i \lambda_k \lambda_r} u_{r \alpha k} u_{\alpha i k}-\sum_k \sum_i \sum_{\alpha \neq i}\frac{u_i^2}{\lambda_k \lambda_\alpha \lambda^2_i}\left( u_{\alpha \alpha k} u_{iik}-u^2_{ki\alpha}\right)\\
&=\sum_k \sum_i \sum_{r \neq i} \sum_{\alpha}\frac{u_i u_r}{\lambda_{\alpha}\lambda_i \lambda_k \lambda_r} u_{r \alpha k} u_{\alpha i k}-\sum_k \sum_i \sum_{\alpha \neq i}\frac{u_i^2}{\lambda_k \lambda_\alpha \lambda^2_i} u_{\alpha \alpha k} u_{iik}+ \sum_k \sum_i \sum_{\alpha \neq i}\frac{u_i^2}{\lambda_k \lambda_\alpha \lambda^2_i}u^2_{ki\alpha}\\
&=\sum_k \sum_i \sum_{r \neq i} \sum_{\alpha}\frac{u_i u_r}{\lambda_{\alpha}\lambda_i \lambda_k \lambda_r} u_{r \alpha k} u_{\alpha i k}+\sum_k \sum_i \sum_{\alpha }\frac{u_i^2}{\lambda_k \lambda_\alpha \lambda^2_i}u^2_{ki\alpha} -\sum_k \sum_i \frac{u_i^2}{\lambda_k \lambda^3_i}u^2_{kii}\\
&\quad -\sum_k \sum_i \frac{u_i^2}{\lambda_k \lambda^2_i}u_{iik} \left(\sum_{\alpha \neq i}\frac{u_{\alpha \alpha k}}{\lambda_{\alpha}} \right),
\end{split}
\end{gather*}
and finally from \eqref{rel} we deduce
\begin{equation*}
\begin{split}
L \varphi &= \sum_k \sum_i \sum_{r \neq i} \sum_{\alpha}\frac{u_i u_r}{\lambda_{\alpha}\lambda_i \lambda_k \lambda_r} u_{r \alpha k} u_{\alpha i k}+\sum_k \sum_i \sum_{\alpha }\frac{u_i^2}{\lambda_k \lambda_\alpha \lambda^2_i}u^2_{ki\alpha} \\
&=\sum_k \sum_i \sum_{r} \sum_{\alpha} \frac{u_i u_r}{\lambda_{\alpha}\lambda_i \lambda_k \lambda_r} u_{r \alpha k} u_{\alpha i k} = \sum_k \sum_{\alpha}\frac{1}{\lambda_{\alpha} \lambda_k}\left(\sum_\gamma \frac{u_{\gamma} u_{\alpha k \gamma}}{\lambda_{\gamma}}\right)^2  \ge 0.
\end{split}
\end{equation*}

\medskip
\noindent {\bf Claim 3: the function $\varphi $ defined in \eqref{fi} is constant in $\Omega$.} 

\noindent By Claim 2 and maximum principle for linear elliptic operators, $\varphi$ attains its maximum over $\bar \Omega$ on the boundary $\partial\Omega$. Therefore either
\begin{itemize}
\item[(i)] $\varphi<c$ in $\Omega$, or 
\item[(ii)]$\varphi\equiv c.$ 
\end{itemize}
Suppose by contradiction that $\varphi$ satisfies (i), that is
\begin{equation}
\label{assurdo}
H_{n-1}|Du|^{n+1}-2u <c \qquad \text{in }\Omega.
\end{equation}
Integrating on $\Omega$ both sides in \eqref{assurdo} we get
\begin{equation}
\label{integro}
\int_{\Omega}\left(H_{n-1}|Du|^{n+1}-2u \right) <c |\Omega|.
\end{equation}
Moreover by \eqref{H-S} and \eqref{divk}  we have
\begin{equation}
\label{int_k}
\int_{\Omega}H_{n-1}|Du|^{n+1}= \int_{\Omega} S^{ij}u_i u_j =n \int_{\Omega}(-u) \det(D^2u)=n\int_{\Omega}(-u).
\end{equation} 
Substituting \eqref{int_k} in \eqref{integro} we get
\begin{equation}
\label{contr}
\int_{\Omega} (-u)  < \frac{c}{n+2}|\Omega|.
\end{equation}
On the other hand \eqref{poz} and \eqref{int_k} imply
\begin{equation*}
-\frac{n}{n+1}\int_{\Omega}(-u)  + \frac{c}{n+1} \int_{\de \Omega} \langle x,\nu\rangle=n \int_{\Omega}(-u),
\end{equation*}
that is
\begin{equation*}
\left(n+\frac{n}{n+1}\right)\int_{\Omega}(-u)  =\frac{c}{n+1} \int_{\de \Omega} \langle x,\nu\rangle.
\end{equation*}
By divergence theorem we finally get
\begin{equation*}
\label{fine}
\int_{\Omega}(-u)  \,dx =\frac{c}{n+2}|\Omega|,
\end{equation*}
that is in contradiction with \eqref{contr}.

\medskip
\noindent {\bf Claim 4:  $\Omega$ is an ellipsoid.}

\noindent Being $\varphi$ constant we have 
\begin{equation}
\label{eq_u}
H_{n-1}|Du|^{n+1}-2u=c \quad \text{in }\bar{\Omega},
\end{equation}
with  $c= 2\max_{\bar{\Omega}} (-u)$.

\noindent Let us consider the following positive, increasing function 
\begin{equation}
\label{g}
g(s)=\frac{n+1}{2n}[c^{n/(n+1)}-(c-2s)^{n/(n+1)}], \quad 0\le s < \frac{c}{2}.
\end{equation}
By \eqref{eq_u}, we get that the function $\psi(x)=g(-u(x))$ satisfies the following equation 
\begin{equation}
\label{eq_g}
H_{n-1}|D\psi|^{n+1}=1 \quad \text{in } \bar\Omega.
\end{equation}
Denote by 
\begin{equation*}
\label{livelli}
\Omega(t)=\{x \in \Omega:\> \psi(x)>t\},\qquad t>0;
\end{equation*}
then $\de \Omega(t)$ is the $t$-level set of the function $\psi(x)$. By \eqref{eq_g} we have
$$
\frac{1}{|D\psi|}=H_{n-1}[\de \Omega(t)]^{\frac{1}{n+1}}.
$$
Thus the family $\Omega(t)$ is a one-parameter family of solutions to the affine curvature flow \eqref{evol}. By Theorem \ref{lemma_Andrews} estimate \eqref{derandr} holds, that is
\begin{equation}\label{bi1}
\frac{\de}{\de t} \left(|\Omega(t)|^{-\frac{n-1}{n+1}} \int_{\de \Omega(t)}H_{n-1}^{\frac{1}{n+1}}\right)\ge 0.
\end{equation}
We claim that equality sign holds in \eqref{bi1}. Hence by  Theorem \ref{lemma_Andrews} $\Omega(t)$ is an ellipsoid for any $0<t<t_E$ and in particular we get that $\Omega(0)=\Omega$ is an ellipsoid. 

Let us define  
\begin{equation*}
\nu(t)=|\Omega(t)|=|\{x \in \Omega \colon \psi(x)>t\}| =|\{x\in\Omega: g(-u(x))>t\}|,\qquad t>0,
\end{equation*}
and
\begin{equation*}
\mu(s)=|\{x \in \Omega \colon -u(x)>s\}| , \qquad s>0;
\end{equation*}
clearly  $\nu(g(s))=\mu(s)$. Taking into account \eqref{eq_g} we immediately have   
\begin{equation}
\label{der_nu}
\nu'(g(s))g'(s)=\mu'(s)=-\int_{\psi=g(s)} \frac{1}{|D\psi|}=-\int_{-u=s}H_{n-1}^{\frac{1}{n+1}}.
\end{equation}
On the other hand from \eqref{eq_u} we deduce 
$$
\mu(s)=\frac{1}{n}\int_{-u>s}\det D^2 u=\frac{1}{n}\int_{-u=s}H_{n-1} |Du|^n=\frac{(c-2s)}{n} \int_{-u=s} \frac{1}{|Du|}=-\mu'(s)\frac{(c-2s)}{n}.
$$
The above equality and \eqref{g}  give
$$
\mu'(s)=-n \mu(s)(g'(s))^{n+1},
$$
and hence
\begin{equation}
\label{espr_der}
\nu'(g(s))=-n \mu(s)g'(s)^{n}=-n \nu(g(s)) \left[\frac{2n}{n+1}g(s)-c^{n/(n+1)}\right]^{-1}.
\end{equation}
Using \eqref{der_nu} and \eqref{espr_der} we finally obtain
\begin{gather*}
\begin{split}
&\frac{\de}{\de g} \left( \nu(g)^{-\frac{n-1}{n+1}}\int_{\psi=g}H_{n-1}^{\frac{1}{n+1}}\right)=-\frac{\de}{\de g} \left( \nu(g)^{-\frac{n-1}{n+1}} \nu'(g)\right)
\\
&=n\frac{\de}{\de g} \left[ \nu(g)^{\frac{2}{n+1}} \left(\frac{2n}{n+1}g-c^{n/(n+1)}\right)^{-1} \right]\\
&=n\left[ \frac{2}{n+1}\nu(g)^{-\frac{n-1}{n+1}}\nu'(g)\left(\frac{2n}{n+1}g-c^{n/(n+1)}\right)^{-1}- \nu(g)^{\frac{2}{n+1}}\frac{2n}{n+1}\left(\frac{2n}{n+1}g-c^{n/(n+1)}\right)^{-2}\right]=0.
\end{split}
\end{gather*}

\begin{remark}\label{rem2}
When $n=2$ it is possible to skip the use of the affine mean curvature flow in the proof of Claim 4. Indeed, being $\varphi$ constant in $\Omega$, from \eqref{eq2}
$$
u_1^2(-\det D^2 u_2)+u_2^2(-\det D^2 u_1)-u_1u_2(u_{111}u_{222}-u_{112}u_{122})=0.
$$
On the other hand by \eqref{alphabeta} $$\alpha^2(-\det D^2 u_2)+\beta^2(-\det D^2 u_1)-\alpha\beta(u_{111}u_{222}-u_{112}u_{122})\ge 0$$
for all $\alpha,\beta\in\R$. Hence we have
$$(u_{111}u_{222}-u_{112}u_{122})^2=4(\det D^2 u_1)(\det D^2 u_2).$$
Then, for almost every $x\in\Omega$ we can use a reference frame where $u_1=0$ and $u_2\ne 0$ and we deduce that $\det D^2 u_1=0$ and $u_{111}u_{222}-u_{112}u_{122}=0$. These identities, together with \eqref{der_det}, enforce all third order derivatives of $u$ to vanish at $x$. By continuity the same holds in the whole $\Omega$, and the claim is proved. 
\end{remark}

\end{document}